\catcode`\^^Z=9
\catcode`\^^M=10
\output={\if N\header\headline={\hfill}\fi
\plainoutput\global\let\header=Y}
\magnification\magstep1
\tolerance = 500
\hsize=14.4true cm
\vsize=22.5true cm
\parindent=6true mm\overfullrule=2pt
\newcount\kapnum \kapnum=0
\newcount\parnum \parnum=0
\newcount\procnum \procnum=0
\newcount\nicknum \nicknum=1
\font\ninett=cmtt9

\font\ninebf=cmbx9

\font\sixbf=cmbx6
\font\ninesl=cmsl9

\font\nineit=cmti9

\font\ninerm=cmr9

\font\sixrm=cmr6
\font\ninei=cmmi9
\font\eighti=cmmi8
\font\sixi=cmmi6
\skewchar\ninei='177 \skewchar\eighti='177 \skewchar\sixi='177
\font\ninesy=cmsy9
\font\eightsy=cmsy8
\font\sixsy=cmsy6
\skewchar\ninesy='60 \skewchar\eightsy='60 \skewchar\sixsy='60
\font\titelfont=cmr10 scaled 1440
\font\paragratit=cmbx10 scaled 1200

\font\name=cmcsc10
\font\emph=cmbxti10

\font\tenmsbm=msbm10
\font\sevenmsbm=msbm7
%

%
\font\got=eufm10

\font\teneufm=eufm10
\font\seveneufm=eufm7
\font\fiveeufm=eufm5
\newfam\eufmfam
\textfont\eufmfam=\teneufm
\scriptfont\eufmfam=\seveneufm
\scriptscriptfont\eufmfam=\fiveeufm

\font\tenmsam=msam10
\font\sevenmsam=msam7
\font\fivemsam=msam5
\newfam\msamfam
\textfont\msamfam=\tenmsam
\scriptfont\msamfam=\sevenmsam
\scriptscriptfont\msamfam=\fivemsam
\font\tenmsbm=msbm10
\font\sevenmsbm=msbm7
\font\fivemsbm=msbm5
\newfam\msbmfam
\textfont\msbmfam=\tenmsbm
\scriptfont\msbmfam=\sevenmsbm
\scriptscriptfont\msbmfam=\fivemsbm
\def\Bbb#1{{\fam\msbmfam\relax#1}}
\def\cz{{\kern0.4pt\Bbb C\kern0.7pt}
}
\def\ez{{\kern0.4pt\Bbb E\kern0.7pt}
}
\def\fz{{\kern0.4pt\Bbb F\kern0.3pt}}
\def\gz{{\kern0.4pt\Bbb Z\kern0.7pt}}
\def\hz{{\kern0.4pt\Bbb H\kern0.7pt}
}
\def\kz{{\kern0.4pt\Bbb K\kern0.7pt}
}
\def\nz{{\kern0.4pt\Bbb N\kern0.7pt}
}
\def\oz{{\kern0.4pt\Bbb O\kern0.7pt}
}
\def\rz{{\kern0.4pt\Bbb R\kern0.7pt}
}
\def\sz{{\kern0.4pt\Bbb S\kern0.7pt}
}
\def\pz{{\kern0.4pt\Bbb P\kern0.7pt}
}
\def\qz{{\kern0.4pt\Bbb Q\kern0.7pt}
}
\newskip\ttglue
\def\ninepoint{\def\rm{\fam0\ninerm}%
  \textfont0=\ninerm \scriptfont0=\sixrm \scriptscriptfont0=\fiverm
  \textfont1=\ninei \scriptfont1=\sixi \scriptscriptfont1=\fivei
  \textfont2=\ninesy \scriptfont2=\sixsy \scriptscriptfont2=\fivesy
  \textfont3=\tenex \scriptfont3=\tenex \scriptscriptfont3=\tenex
  \def\it{\fam\itfam\nineit}%
  \textfont\itfam=\nineit
  \def\sl{\fam\slfam\ninesl}%
  \textfont\slfam=\ninesl
  \def\bf{\fam\bffam\ninebf}%
  \textfont\bffam=\ninebf \scriptfont\bffam=\sixbf
   \scriptscriptfont\bffam=\fivebf
  \def\tt{\fam\ttfam\ninett}%
  \textfont\ttfam=\ninett
  \tt \ttglue=.5em plus.25em minus.15em
  \normalbaselineskip=11pt
  \font\name=cmcsc9
  \let\sc=\sevenrm
  \let\big=\ninebig
  \setbox\strutbox=\hbox{\vrule height8pt depth3pt width0pt}%
  \normalbaselines\rm
  \def\sl{\it}}

\headline={\ifodd\pageno\rightheadline\else\leftheadline\fi}
\def\rightheadline{\ninepoint Paragraphen"uberschrift\hfill\folio}
\def\leftheadline{\ninepoint\folio\hfill Chapter"uberschrift}
\let\header=Y
\def\titel#1{\need 9cm \vskip 2truecm
\parnum=0\global\advance \kapnum by 1
{\baselineskip=16pt\lineskip=16pt\rightskip0pt
plus4em\spaceskip.3333em\xspaceskip.5em\pretolerance=10000\noindent
\titelfont Chapter \uppercase\expandafter{\romannumeral\kapnum}.
#1\vskip2true cm}\def\leftheadline{\ninepoint
\folio\hfill Chapter \uppercase\expandafter{\romannumeral\kapnum}.
#1}\let\header=N
}
\def\Titel#1{\need 9cm \vskip 2truecm
\global\advance \kapnum by 1
{\baselineskip=16pt\lineskip=16pt\rightskip0pt
plus4em\spaceskip.3333em\xspaceskip.5em\pretolerance=10000\noindent
\titelfont\uppercase\expandafter{\romannumeral\kapnum}.
#1\vskip2true cm}\def\leftheadline{\ninepoint
\folio\hfill\uppercase\expandafter{\romannumeral\kapnum}.
#1}\let\header=N
}
\def\need#1cm {\par\dimen0=\pagetotal\ifdim\dimen0<\vsize
\global\advance\dimen0by#1 true cm
\ifdim\dimen0>\vsize\vfil\eject\noindent\fi\fi}
\def\neupara#1{\par\penalty-2000
\procnum=0\global\advance\parnum by 1
\vskip1cm\noindent{\paragratit \the\parnum. #1}%
\def\rightheadline{\ninepoint\S\the\parnum.\ #1\hfill \folio}%
\vskip 8mm\noindent}
\def\Proclaim #1 #2\finishproclaim {\bigbreak\noindent
{\bf#1\unskip{}. }{\it#2}\medbreak\noindent}
%
\gdef\proclaim #1 #2 #3\finishproclaim {\bigbreak\noindent%
\global\advance\procnum by 1
{%
{\relax\ifodd \nicknum
\hbox to 0pt{\vrule depth 0pt height0pt width\hsize
   \quad \ninett#3\hss}\else {}\fi}%
\bf\the\parnum.\the\procnum\ #1\unskip{}. }
{\it#2}
\immediate\write\num{\string\def
 \expandafter\string\csname#3\endcsname
 {\the\parnum.\the\procnum}}
\medbreak\noindent}
\newcount\stunde \newcount\minute \newcount\hilfsvar
\def\uhrzeit{
    \stunde=\the\time \divide \stunde by 60
    \minute=\the\time
    \hilfsvar=\stunde \multiply \hilfsvar by 60
    \advance \minute by -\hilfsvar
    \ifnum\the\stunde<10
    \ifnum\the\minute<10
    0\the\stunde:0\the\minute~Uhr
    \else
    0\the\stunde:\the\minute~Uhr
    \fi
    \else
    \ifnum\the\minute<10
    \the\stunde:0\the\minute~Uhr
    \else
    \the\stunde:\the\minute~Uhr
    \fi
    \fi
    }
\def\calA{{\cal A}} \def\calB{{\cal B}}

\def\calE{{\cal E}} 
\def\calG{{\cal G}} \def\calH{{\cal H}}
\def\calI{{\cal I}}

\def\calO{{\cal O}} \def\calP{{\cal P}}
 
\def\calS{{\cal S}} \def\calT{{\cal T}}
\def\calU{{\cal U}} 
 \def\calX{{\cal X}}
 \def\calZ{{\cal Z}}

\def\gotm{\hbox{\got m}}

\def\dim{\mathop{\rm dim}\nolimits}

\def\GL{\mathop{\rm GL}\nolimits}

\def\mod{\mathop{\rm mod}\nolimits}
\def\O{{\rm O}}

\def\SL{\mathop{\rm SL}\nolimits}

\def\Sp{\mathop{\rm Sp}\nolimits}

\def\Sym{\mathop{\rm Sym}\nolimits}
\def\boxit#1{
  \vbox{\hrule\hbox{\vrule\kern6pt
  \vbox{\kern8pt#1\kern8pt}\kern6pt\vrule}\hrule}}
\def\Boxit#1{
  \vbox{\hrule\hbox{\vrule\kern2pt
  \vbox{\kern2pt#1\kern2pt}\kern2pt\vrule}\hrule}}

\def\zwischen#1{\bigbreak\noindent{\bf#1\medbreak\noindent}}

\def\smallni{\smallskip\noindent }
\def\medni{\medskip\noindent }
\def\bigni{\bigskip\noindent }
\def\Isom{\mathop{\;{\buildrel \sim\over\longrightarrow }\;}}
\def\lo{\longrightarrow}

\def\loma{\longmapsto}
\def\betr#1{\vert#1\vert}

\def\imag{{\rm i}}
\def\pii{\pi {\rm i}}
\def\veps{\varepsilon}
\def\square{\hbox{\hbox to 0pt{$\sqcup$\hss}\hbox{$\sqcap$}}}
\def\qed{\ifmmode\square\else{\unskip\nobreak\hfil
\penalty50\hskip3em\null\nobreak\hfil\square
\parfillskip=0pt\finalhyphendemerits=0\endgraf}\fi}
\def\pn{\the\parnum.\the\procnum}
\def\downmapsto{{\buildrel
        {\vbox{\hbox{\hskip.2pt$\scriptstyle-$}}}
        \over{\raise7pt\vbox{\vskip-4pt\hbox{$\textstyle\downarrow$}}}}}

\nopagenumbers
\immediate\newwrite\num
\nicknum=0  
\let\header=N
\def\transpose#1{\kern1pt{^t\kern-1pt#1}}%

\def\Cl{{\rm Cl}}
\def\Clan{\Cl_{\hbox{\sevenrm an}}}
\def\Qan{Q_{\hbox{\sevenrm an}}}
\def\Hilb{\hbox{\rm Hilb}}
\def\RAND#1{\hbox to 0mm{\hss\vtop to 0pt{%
  \raggedright\ninepoint\parindent=0pt%
  \baselineskip=1pt\hsize=2cm #1\vss}}\noindent}
  \def\NumOrb{1.1}

\def\PresRu{1.3}
\def\RulStan{1.4}
\def\KoQua{2.1}
\def\ClNode{2.2}
\def\SnCus{3.1}
\def\NodeTrans{3.2}
\def\KonvB{4.1}
\def\MultBorch{4.2}

\def\CocB{4.5}

\def\UgeN{5.1}

\def\RelCl{5.6}

\def\nPPn{6.1}
\def\HauptS{6.2}

\def\freiCM{7.2}
\def\FixSet{7.3}

\def\ordZw{8.2}
 
\noindent
\centerline{\titelfont On Siegel threefolds with a projective Calabi--Yau model}%
\def\leftheadline{\ninepoint\folio\hfill
On Siegel three folds with a projective Calabi--Yau model}%
\def\rightheadline{\ninepoint Introduction\hfill \folio}%
\headline={\ifodd\pageno\rightheadline\else\leftheadline\fi}
\vskip 1.5cm
\leftline{\it \hbox to 6cm{Eberhard Freitag\hss}
Riccardo Salvati
Manni  }
  \leftline {\it  \hbox to 6cm{Mathematisches Institut\hss}
Dipartimento di Matematica, }
\leftline {\it  \hbox to 6cm{Im Neuenheimer Feld 288\hss}
Piazzale Aldo Moro, 2}
\leftline {\it  \hbox to 6cm{D69120 Heidelberg\hss}
 I--00185 Roma, Italy. }
\leftline {\tt \hbox to 6cm{freitag@mathi.uni-heidelberg.de\hss}
salvati@mat.uniroma1.it}
\vskip1cm
\centerline{\paragratit \rm  2011}%
\vskip5mm\noindent%
\let\header=N%
\def\imag{{\rm i}}%
{\paragratit Introduction}%
\medni
In the papers [FS1], [FS2]
we described some Siegel modular threefolds which admit
a weak Calabi--Yau model.\footnote{*)}{\ninepoint
In this paper a weak Calabi--Yau threefold is understood as
a compact complex threefold of whose first Betti number vanishes and which admits
an everywhere holomorphic differential form of degree 3 without zeros.}
Not all of them admit a {\it projective\/} model.
In fact, Bert van Geemen, in a private communication,
pointed out a significative example which cannot admit
a projective model.   His comment was a starting motivation
for this paper.
We mention that a weak Calabi--Yau
threefold is projective if, and only if, it admits a Kaehler metric.
The purpose of this paper is to exhibit  criterions for the projectivity,
to treat several
examples, and
to compute their Hodge numbers.
\smallskip
Basic for our examples is a certain complete intersection $\calX$
of four quadrics introduced
the paper [GN] of van Geemen and Nygaard:
$$\eqalign{Y_0^2&=X_0^2+X_1^2+X_2^2+X_3^2,\cr
Y_1^2&=X_0^2-X_1^2+X_2^2-X_3^2,\cr
Y_2^2&=X_0^2+X_1^2-X_2^2-X_3^2,\cr
Y_3^2&=X_0^2-X_1^2-X_2^2+X_3^2\cr}\leqno\qquad\calX:$$
The variety $\calX$ has 96 isolated singularities which are ordinary
double points (nodes). One of them, called the {\it standard node,\/} is
$$\eta=[\sqrt 2,0,\sqrt 2,0,\,1,1,0,0].$$
In the paper [CM] it has
been pointed out that the results of [GN] imply that
$\calX$ admits a resolution that is a (projective)
Calabi--Yau threefold. The holomorphic three form
without zeros (unique up to a constant) on this model is given by
$${X_2^4\over Y_0Y_1Y_2Y_3}d(X_0/X_2)
\wedge d(X_1/X_2)\wedge d(X_3/X_2).$$
The basic result -- essentially due to van Geemen and
Nygaard [GN] -- is the following theorem.
\Proclaim
{Theorem}
{The Hodge numbers of a Calabi--Yau desingularization of $\calX$ are
$$h_{11}=32,\qquad h_{12}=0.$$
Hence this Calabi--Yau manifold is rigid.}
\finishproclaim
We recall that the two Hodge numbers and the Euler number $e$ are
related by the formula
$$e=2(h_{11}-h_{12}),$$
hence the Euler number for this example is $e=64$.
\smallskip
In [FS2] we introduced a finite group $\calG\subset\GL(8,\cz)$
which acts on $\calX$. It is generated by the four transformations
$$\eqalign{
U_1:&\qquad\qquad
(Y_0,Y_1,Y_3,Y_2,X_0,X_3,X_2,X_1),\cr
U_2:&\qquad\qquad
(Y_0,Y_3,Y_2,Y_1,X_0,X_1,X_3,X_2),\cr
T:&\qquad\qquad
(Y_0,-\imag Y_1,Y_2,-\imag Y_3,X_1, X_0,X_3,X_2),\cr
J:&\qquad\qquad
\hbox to 0pt{\hss$\sqrt2\cdot$}(X_0,X_1,X_2,X_3,
Y_0/2,Y_1/2,Y_2/2,Y_3/2).\cr
}$$
This group contains the subgroup $\calZ$
of order 4 generated by multiplication with
$\imag$. The group
$$\bar\calG:=\calG/\calZ$$
acts faithfully on $\calX$.
Its order is $24\,576=2^{13}\cdot 3$.
There is a subgroup $\calH\subset\calG$ of index two which leaves the
holomorphic three form invariant, namely, the group generated by
$$U_1U_2,\quad U_1T,\quad U_2T,\quad J.$$
It also contains $\calZ$ and
$$\bar\calH:=\calH/\calZ$$
is a group of order
$12\,288=2^{12}\cdot 3$ which acts faithfully on $\calX$ and which leaves the
holomorphic 3-form invariant.
In [FS2] has been proved:
\Proclaim
{Theorem}
{Let $G$ be any subgroup of $\bar\calH$. Then there exists a
desingularization (in the category of complex spaces)
of $\calX/G$ that is a weak Calabi--Yau
threefold.}
\finishproclaim
Actually we proved more in [FS2], namely that the varieties
$\calX/G$ are Siegel threefolds.
\smallskip
There are $4\,117$ conjugacy classes of subgroups of $\bar\calH$
so, in principle, we get $4\,117$ examples of weak Calabi--Yau threefolds.
Of course there might be biholomophic equivalent ones under them, but, this may be
difficult to decide.
\smallskip
Let $G\subset\bar\calH$ be a subgroup. One can ask whether there exists a
resolution of $\calX/G$ in the form of a {\it projective\/} Calabi--Yau manifold.
In all cases that we are able to treat we will construct first a
resolution
$\tilde\calX\to\calX$ such that $G$ extends to a group of biholomorphic
self maps of $\tilde\calX$. After that we will construct a resolution
of $\tilde\calX/G$ in the form of a weak Calabi--Yau manifold.
\smallskip
We explain the construction of $\tilde\calX$. For this we have to consider
two classes $\calA$, $\calB$ of nodes of $\calX$.
\smallskip
A node $a$ belong to the class $\calA$
if its stabilizer $G_a$ contains an element $g$ that is conjugate to the
transformation
$$(Y_0,\dots, X_3)\loma (Y_0,-Y_1,Y_2,-Y_3,X_0,X_1,-X_2,-X_3).$$
This transformation occurred already in [FS2]. It has the following
property. Consider the blow up of the node. The exceptional divisor
($\cong \pz^1\times \pz^1$) is in the fixpoint locus of $g$.
\smallskip
Next we explain the second class $\calB$. For this one needs the small resolutions
of the node. A resolution is called {\it small\/} if the exceptional set
is a curve. One knows that a node admits two (isomorphy classes)
of small resolutions.
By a {\it ruling\/} of the node we understand the choice of one of
the two small resolutions.
There is a subgroup of index two
of $\bar\calG_a$ that preserves the rulings.
We described this subgroup in [FS2].
This subgroup extends to the small resolutions of the node.
By definition, the class $\calB$ consists of all nodes $a$ such that
the elements of the stabilizer $G_a$ preserve the rulings of the node.
Assume that $\calA\cup\calB$ is the set of all nodes.
Then we can construct a resolution $\tilde\calX\to\calX$ as follows.
For the nodes in $\calA$ we take the blow up and for the rest we take
a small resolution. This can be done in such a way that $G$ acts on
$\tilde\calX$ as group of biholomorphic transformations.
Now a general theorem, essentially due to Roan [Ro], shows the following
result.
(compare [FS2], Theorem 1.5).
\Proclaim
{Theorem}
{Let $G\subset\bar\calH$ be a subgroup such that each
node is contained in $\calA\cup\calB$. Then their exists
a resolution $\tilde\calX\to\calX$ in the category of complex spaces
such that the action of $G$
extends to $\tilde\calX$ and such that $\tilde\calX/G$ admits a
resolution in the form of a weak Calabi--Yau threefold.
There are\/ $2\,791$ conjugacy classes of
subgroups of $\bar\calH$ with this property.}
\finishproclaim
Now there arises the question whether $\tilde\calX$ can be obtained as
{\it projective variety.\/} If this is the case then we can construct
the resolution $\tilde\calX/G$
on the form of the Hilbert scheme $G$-$\Hilb(\tilde\calX)$
as a (projective) Calabi--Yau manifold (see [BKR],
compare also [FS1],  Theorem 2.6).  We will derive a projectivity
criterion to obtain the following result.
\Proclaim
{Theorem}
{Let $G\subset\bar\calH$ be a group such that $\calA\cup\calB$ is the set
of all nodes. Assume that for each node $a\in\calB$, $a\not\in\calA$,
there exists a divisor on $\calX$ that runs in a non trivial way into the node.
Then $\tilde\calX$ can be constructed as projective manifold. As a consequence,
there exists a resolution of $\calX/G$ in the form of a (projective)
Calabi-Yau manifold. There are $1\,986$ conjugacy classes of
subgroups of $\bar\calH$ with this property.}
\finishproclaim
Here ``non-tivial'' means that the divisor is not the divisor of a meromorphic
function in any small (analytic) open neighborhood of the node.
To prove projectivity, one needs information about the map $\Cl(\calX)\to\Cl(\calX,a)$.
Here $\Cl(\calX)$ is the group of divisor classes of $\calX$, and
the group $\Cl(\calX,a)$ is the limit of the analytic divisor class groups
of small neighborhoods of $a$.
We already have seen
in [FS2] that $\Cl(\calX)$ is a group of rank 32. In Sect.~2 we describe an explicit system of
generators.  The group $\Cl(\calX,a)$ is isomorphic to $\gz$ for a node.
A projective small resolution of a node can be obtained if one blows up a divisor
that runs non-trivially into the node which means that the image in
$\Cl(\calX,a)$ is different from 0. Hence we have to describe the images of the
generating divisor in $\Cl(\calX,a)$. This is done in Sects.~3--5 and needs analytic
methods. Actually we use the description
of van Geemen and Nygaard as modular variety.
The nodes correspond to certain 0-dimensional cusps. We develop a theory
of local Borcherds products to describe $\Cl(\calX,a)$. This theory enables us to
compute the map $\Cl(\calX)\to\Cl(\calX,a)$ explicitly.
It looks rather involved to introduce a new theory of local Borcherds products
in this context. But this theory may be of interest in its own right.
\smallskip
In the last three sections we show how the Hodge numbers of many examples
can be computed.
Many examples rest on computer calculations. Nevertheless we treated some
examples in detail where the calculations can be done by hand.
\neupara{The divisor class group}%
In the paper [FS2] we investigated  the divisor class group $\Cl(\calX)$
and proved that $\Cl(\calX)\otimes_\gz\qz$ has dimension 32.
Generators have been found using the decomposition of the divisor
of Igusa's modular form $\chi_{35}$ into irreducible components.
We reproduce these results in a modified form.
For this we start with the three forms
$$X_2-X_3,\quad X_0X_2+X_1X_3,\quad  X_0-X_1-X_2-X_3.$$
Their zero divisors on $\calX$ are not irreducible.
Each of them can be decomposed into two divisors as follows:
$$\eqalign{
D_1^\pm:\quad &X_2-X_3=Y_1\pm Y_3=0,\cr
D_2^\pm:\quad &X_0X_2+X_1X_3=Y_0Y_1\pm Y_2Y_3=0,\cr
D_3^\pm:\quad &X_0-X_1-X_2-X_3=Y_1(X_1+X_3)\pm(\sqrt2/2) Y_2Y_3.\cr
}$$
To be precise we mention that this is only a set theoretical description.
The precise definition of -- for example $D_3^+$ -- is that the associated
ideal the divisor is the radical of the image of the ideal
$$(X_0-X_1-X_2-X_3,Y_1(X_1+X_3)+(\sqrt2/2) Y_2Y_3)$$
in the coordinate ring of $\calX$ (factor ring of $\cz[Y_0,\dots,X_3]$ by the
defining ideal). We will make use of the orbits of the divisors $D_i^\pm$
under $\calG$.
\proclaim
{Lemma}
{For every $i\in\{1,2,3\}$ the divisor $D_i^+$ and $D_i^-$ are in the same
$\calG$-orbit.
The order of the orbit of the divisors $D_i^\pm$ under $\calG$ is
$48$ for $D_1^\pm$, $12$ for $D_2^\pm$, and $128$ for $D_3^\pm$.}
NumOrb%
\finishproclaim
From [FS2] we recall the result that $\Cl(\calX)\otimes_\gz\qz$ is the direct
sum of $\calG$-irreducible subspaces of dimension $1$, $3$, $12$, $16$.
We consider the subspaces of $\Cl(\calX)\otimes_\gz\qz$ generated by the orbits
of $D_i$. Each of this three subspace contains the trivial one-dimensional
representation.
The results described in [FS2] imply the following proposition.
\proclaim
{Proposition}
{Consider the factor space of $\Cl(\calX)\otimes_\gz\qz$ by the one-dimensional
representation. Its three irreducible components can be described as
the $\calG$-submodules generated by the divisors $D_i^\pm$. More precisely,
$D_1^\pm$ generates the subspace of dimension $12$,
$D_2^\pm$ generates the subspace of dimension
$3$, and $D_3^\pm$ generates the subspace of dimension $16$.}
IrrComp%
\finishproclaim
\indent
We want to describe how the generating divisors run into the standard node.
For this we need information about its stabilizer.
It has been described in [FS2].
\proclaim
{Lemma}
{The following transformations stabilize the standard node
and preserves its ruling:
\smallni
\halign{\hskip1cm$#$\hfil\cr
(Y_0, Y_3, Y_2, Y_1, X_0, X_1, X_3, X_2),\cr
(Y_0, \imag Y_1, Y_2, \imag Y_3, X_1, X_0, X_3, X_2),\cr
(Y_2, Y_3, Y_0, Y_1, X_0, X_1, -\imag X_2, -\imag X_3),\cr
(Y_0, Y_1, Y_2, -Y_3, X_0, X_1, X_2, -X_3),\cr
(Y_0, -Y_1, Y_2, -Y_3, X_0, X_1, X_2, X_3),\cr
(Y_0, Y_1, Y_2, Y_3, X_0, X_1, -X_2, -X_3),\cr
(2X_0, 2X_2, 2X_1, 2X_3, Y_0, Y_2, Y_1, Y_3).\cr}
\smallni
Their images in $\bar\calG$  generate a subgroup order $128=2^7$
which is a subgroup of index two of
the stabilizer.
}
PresRu%
\finishproclaim
Now we consider the set of all basic divisors that run into the standard node
and decompose it into orbits.
\proclaim
{Theorem}
{Consider the set of $188$ divisors described in Lemma \NumOrb.
Let $\calS$ be the subset of divisors that run into the standard node.
The group of all elements of $\calG$ that stabilize the standard node and its
rulings (see \PresRu) acts on this set. There are six orbits that can be represented
be the six divisors
$$D_1^+,D_1^-,D_2^+,D_2^-,D_3^+,D_3^-.$$}
RulStan%
\finishproclaim
\vskip-\bigskipamount\noindent
In the next sections we will use this result to determine the classes
of the basic divisors in small neighborhood of a node.
\neupara{The local divisor class group of a node}%
Let $M$ be a locally compact space. An open subset $U$ is called a neighborhood
of $\infty$ if its complement is a compact subset.
We define
$$H^q_\infty(M,\gz):=\lim_{\longleftarrow}H^q(U,\gz)$$
as the limit of the cohomology groups of all open neighborhoods of $\infty$.
There is an exact sequence
$$\cdots\lo H^q_c(M,\gz)\lo H^q(M,\gz)\lo H^q_\infty(M,\gz)\lo H^{q+1}_c(M,\gz)\lo\cdots$$
We mention that it is sufficient that $U$ runs through a
fundamental system of open neighborhoods of $\infty$.
\smallskip
We apply this to the tangent bundle $M$ of the sphere $S_n$ for $n>1$. Explicitly, this is
$$M=\Big\{(x,y)\in\rz^{2(n+1)};\quad \sum x_i^2=1,\ \sum x_iy_i=0\Big\}.$$
Since this is homotopically equivalent to the sphere, we have
$$H^q(M,\gz)=\cases{\gz &if $q=0$ or $q=n$,\cr 0& else.}$$
Poincar\`e duality gives (for the manifold $M$ of dimension  $2n$)
$$H^q_c(M,\gz)=\cases{\gz &if $q=n$ or $q=2n$,\cr 0& else.}$$
If $H^q(M,\gz)$ is different from 0, its generator has not compact support.
Hence $H^q_c(M,\gz)\to H^q(M,\gz)$ is always the zero map.
Now the exact sequence above shows the following.
$$H^q_\infty(M,\gz)=\cases{\gz &if $q\in\{0,n-1,n,2n-1\}$,\cr 0& else.}$$
We consider for a positive $C$
$$\eqalign{
M(C)&=\{(x,y)\in M;\quad \sum y_i^2=C\},\cr
M(>C)&=\{(x,y)\in M;\quad \sum y_i^2>C\}.\cr}$$
All these sets are homotopically equivalent. Moreover, the sets $M(>C)$
define a fundamental system of open neighborhoods of $\infty$.
This shows
$$H^q(M(C),\gz)=H^q_\infty(M,\gz).$$
There is another way to read this result.
We consider the quadric
$$Q=\Big\{z\in\cz^{n+1},\quad \sum z_i^2=0\Big\}.$$
We want to determine the cohomology of $Q-\{0\}$.
This is homotopically equivalent to
$$Q(C)=\Big\{z\in Q,\quad \sum \betr{z_i}^2=C\Big\}.$$
There is a topological map
$$M(1)\lo Q(2),\quad (x,y)\loma z:=x+\imag y.$$
We obtain the following result.
\proclaim
{Lemma}
{Let $Q\subset\cz^{n+1}$ be the quadric defined by
$\sum z_i^2=0$. Then
$$H^q(Q-\{0\},\gz)=\cases{\gz &if $q\in\{0,n-1,n,2n-1\}$,\cr 0& else.}$$}
KoQua%
\finishproclaim
We need an analytic version
$\Clan(X)$ of the divisor class group
for an irreducible normal complex space $X$.
An analytic divisor is a formal linear
combination $D=\sum_Y n_YY$, $n_Y\in\gz$,
of irreducible closed complex subspaces of codimension 1 such that
for each compact subset $K\subset X$ there exist only finitely many $Y$
with $Y\cap K\ne 0$ and $n_Y\ne 0$. Since $X$ is normal, each non-zero meromorphic
function $f$ defines a divisor $(f)$, called a principal
divisor. The divisor class group is the factor group of all
divisors and the subgroup of all principal divisors.
If $S\subset X$ is a closed analytic subspace of codimension $\ge 2$ then
$\Clan(X)=\Clan(X-S)$. If $X$ is a projective variety then $\Clan(X)$
and the algebraic divisor class group $\Cl(X)$ are naturally isomorphic.
\smallskip
We need a local variant of the analytic divisor class group.
Let $X$ be an  normal complex space and let $a\in X$. We set
$$\Clan(X,a):=\lim_U\Clan(U),$$
where $U$ runs through all
connected open neighborhoods $U$ of $a$.
It is of course sufficient that $U$ runs through a fundamental system of neighborhoods.
\proclaim
{Lemma}
{Let
$$Q:=\{z\in\cz^4;\quad z_1z_4=z_2z_3\}.$$
We have
$$\Clan(Q,0)\cong\gz.$$
A generator can be given by the divisor
which is defined by $z_1=z_2=0$. Its negative can be defined by
$z_1=z_3=0$. The blow ups of these two divisors give the two small resolutions.
}
ClNode%
\finishproclaim
{\it Proof.\/} Since $Q$ is a Cohen-Macaulay variety, we have that the cohomology
$H^q_{\{0\}}(\Qan,\calO_{\Qan})$ with support in  the origin vanishes for $q<3$.
This implies
$$H^1(\Qan-\{0\},\calO)=0.$$
Here $\calO$ denotes the sheaf of analytic functions on $\Qan-\{0\}$.
From the exponential sequence we get
$$0\lo H^1(\Qan-\{0\},\calO^*)\lo H^2(\Qan-\{0\},\gz).$$
From Lemma \KoQua\ we know that  $H^2(\Qan-\{0\})\cong\gz$.
One can check that the line bundle related to the divisor $z_1=z_2=0$ goes
to a generator. This means that we get an isomorphism
$$H^1(\Qan-\{0\},\calO^*)\Isom\gz.$$
There is natural injective map
$$\Clan(Q-\{0\})\lo H^1(\Qan-\{0\},\calO^*).$$
Since $H^1(\Qan-\{0\},\calO^*)$ is generated by the image of the divisor
$z_1=z_2=0$ we obtain that it is an isomorphism.
We can repeat the whole consideration for
$$\Big\{z\in Q;\quad \sum \betr{z_i}^2<\veps\Big\}$$
instead of $Q$. This finishes the proof of Lemma \ClNode.\qed
\neupara{The modular approach}%
It is necessary for us to understand the map
$\Cl(\calX)\to\Cl(\calX,a)$ for the nodes of $\calX$.
For this we have to use the realization of $\calX$ as Siegel threefold and
the nodes as certain 0-dimensional cusps.
\smallskip
Following van Geemen and Nygaard,
we described in [FS2] a certain congruence subgroup $\Gamma'\subset\Sp(2,\gz)$
such that $\calX$ is biholomorphic equivalent to the Satake compactification
$\overline{\hz_2/\Gamma'}$ of $\hz_2/\Gamma'$,
namely
$$\Gamma'=\{M\in \Gamma_2[2,4]\cap\Gamma_{2,0,\vartheta}[4];
\quad \det D\equiv\pm1\;\mod\;8\}.$$
For the notations we refer to [FS2].
\smallskip
The biholomorphic map $\calX\cong\overline{\hz_2/\Gamma'}$
is given by the map that assigns
the variables $Y_0,\dots, Y_3,X_0\dots,X_3$ in this ordering
to the theta functions
\medni
\halign{\qquad$\displaystyle#$\quad\hfil&$\displaystyle#$\quad\hfil
&$\displaystyle#$\quad\hfil&$\displaystyle#$\quad\hfil\cr
\vartheta\Bigl[{00\atop 00}\Bigr](Z),& \vartheta\Bigl[{00\atop 10}\Bigr](Z),&
\vartheta\Bigl[{00\atop 01}\Bigr](Z),& \vartheta\Bigl[{00\atop 11}\Bigr](Z),\cr
\noalign{\vskip3mm}
\vartheta\Bigl[{00\atop 00}\Bigr](2Z),& \vartheta\Bigl[{10\atop 00}\Bigr](2Z),&
\vartheta\Bigl[{01\atop 00}\Bigr](2Z),& \vartheta\Bigl[{11\atop 00}\Bigr](2Z).\cr}
\medni
The nodes correspond to certain zero dimensional boundary component.
\proclaim
{Lemma}
{The standard node corresponds to the boundary point
$$\pmatrix{0&0\cr0&\imag\infty}:=\lim_{t\to+\infty}
\pmatrix{\imag/t&0\cr0&\imag t}.$$}
SnCus%
\finishproclaim
{\it Proof.\/} We can use the embedded involution
$$N:=\pmatrix{0&0&1&0\cr 0&1&0&0\cr-1&0&0&0\cr0&0&0&1}$$
to transform this boundary point to the standard boundary point
$$\pmatrix{\imag\infty&0\cr0&\imag\infty}:=\lim_{t\to+\infty}
\pmatrix{\imag t&0\cr0&\imag t}.$$
We have to transform the modular forms
\medni
\vbox{
\halign{$\displaystyle#$\ \hfil&$\displaystyle#$\ \hfil
&$\displaystyle#$\ \hfil&$\displaystyle#$\ \hfil\cr
Y_0=\vartheta\Bigl[{00\atop 00}\Bigr](Z),& Y_1=\vartheta\Bigl[{00\atop 10}\Bigr](Z),&
Y_2=\vartheta\Bigl[{00\atop 01}\Bigr](Z),& Y_3=\vartheta\Bigl[{00\atop 11}\Bigr](Z),\cr
\noalign{\vskip3mm}
X_0=\vartheta\Bigl[{00\atop 00}\Bigr](2Z),& X_1=\vartheta\Bigl[{10\atop 00}\Bigr](2Z),&
X_2=\vartheta\Bigl[{01\atop 00}\Bigr](2Z),& X_3=\vartheta\Bigl[{11\atop 00}\Bigr](2Z)\cr}
}
\medni
by means of the transformation
$$f\vert N(Z)=\det(CZ+D)^{-1}f(NZ).$$
Standard theta transformation formulas show that the transformed
forms up to a joint constant factor are
$$Y_0':=\vartheta\Bigl[{00\atop 00}\Bigr](Z),\quad
Y_1':=\vartheta\Bigl[{10\atop 00}\Bigr](Z),\quad
Y_2':=\vartheta\Bigl[{00\atop 01}\Bigr](Z),\quad
Y_3':=\vartheta\Bigl[{10\atop 01}\Bigr](Z)$$
and
$$\eqalign{
X_0'={1\over\sqrt2}\bigl(\vartheta\Bigl[{00\atop 00}\Bigr](2Z)+
\vartheta\Bigl[{10\atop 00}\Bigr](2Z)\bigr),\quad
&X_1'={1\over\sqrt2}\bigl(\vartheta\Bigl[{00\atop 00}\Bigr](2Z)-
\vartheta\Bigl[{10\atop 00}\Bigr](2Z)\bigr)\cr
X_2'={1\over\sqrt2}\bigl(\vartheta\Bigl[{01\atop 00}\Bigr](2Z)+
\vartheta\Bigl[{11\atop 00}\Bigr](2Z)\bigr),\quad
&X_3'={1\over\sqrt2}\bigl(\vartheta\Bigl[{01\atop 00}\Bigr](2Z)-
\vartheta\Bigl[{11\atop 00}\Bigr](2Z)\bigr).
\cr}$$
The projective coordinates of the cusp ${0\;\;0\;\choose0\,\imag\infty}$
with respect to
the variables $Y_0,\dots,X_3$ are the same as the coordinates of the transformed
cusp ${\imag\infty\;0\choose\;0\;\;\imag\infty}$ with respect to the new
variables $Y_0',\dots,X_3'$. Using the formula
$$\lim_{t\to\infty}\vartheta\Bigl[{a_1a_2\atop b_1b_2}\Bigr]
\pmatrix{\imag t&0\cr0&\imag t}=\cases{1& if $a_1=a_2=0$,\cr 0&else}$$
we get for the value
$[1,0,1,0,1,\;1/\sqrt2,1/\sqrt2,0,0]$. This is the standard node as claimed in
\SnCus.\qed\smallskip
We denote the conjugated group of $\Gamma'$ by
$$\Gamma'':=N\Gamma'N^{-1}.$$
The transformed forms $Y_0',\dots,X_3'$ are generators of the ring of modular
forms of $\Gamma''$. Of course they satisfy the same relations as the
$Y_0,\dots,X_3$.
One gets
$$\eqalign{\Gamma''=\big\{M\in\Gamma_2[2,4];\quad
&a_{12} \equiv d_{21}\equiv \, 0 \,\mod\, 4,\quad
b_{11}\equiv c_{22} \equiv \,0\, \mod\, 8,\cr
&a_{11}d_{22}- b_{12}c_{21}\equiv \pm 1\, \mod\, 8\big\}.\cr}$$
The Siegel-parabolic subgroup $\calP$ of $\Gamma''$ (defined by $C=0$) consists
of all
$$\pmatrix{E&T\cr 0&E}\pmatrix{\transpose U&0\cr 0&U^{-1}},
\quad U\in\calU,\ T\in\calT.$$
Here $\calT$ denotes the set of all integral matrices
$T=\pmatrix{t_0&t_1\cr t_1&t_2}$  that satisfy
$$t_0\equiv 0\mod 8,\quad t_1\equiv 0\mod 2,\quad t_2\equiv 0\mod 4,
$$
and  $\calU$ is the subgroup of $\GL(2,\gz)$ defined by the congruences
$$b\equiv 0\mod 2,\quad c\equiv 0\mod 4,\qquad U=\pmatrix{a&b\cr c&d}.$$
This group contains non-trivial elements of finite order, for example
the diagonal matrix with entries $1,-1$.
\smallskip
So we have proved the following result.
\proclaim
{Lemma}
{There is a biholomorphic map between $\calX$ and
$\overline{\hz_2/\Gamma''}$. It is defined through the correspondence
$X_i\leftrightarrow X_i'$, $Y_i\leftrightarrow Y_i'$. The standard node
corresponds the standard cusp ${\imag\infty\;0\choose 0\;\imag\infty}$.
The corresponding parabolic group is the  group $\calP$ described above.}
NodeTrans%
\finishproclaim
\neupara{0-dimensional cusps}%
We consider the Siegel parabolic group that consists
of all integral symplectic matrices of genus two
$$M=\pmatrix{A&B\cr C&D},\qquad C=0.$$
Let $\calP$ be subgroup of finite index. For simplicity we assume
that  $\calP$ splits. This means the following.
If $M$ is in $\calP$ then the matrix
$$M=\pmatrix{A&0\cr 0&D}$$
is in $\calP$ too. Let $\calU\subset\GL(n,\gz)$ be the subgroup
of all $U$ such that
$$\pmatrix{\transpose U&0\cr 0&U^{-1}}\in\calP,$$
and let $\calT$ be the set of all integral symmetric matrices
$T$ such that
$$\pmatrix{E&T\cr 0&E}\in\calP.$$
Then $\calP$ consists of all
$$\pmatrix{E&T\cr 0&E}\pmatrix{\transpose U&0\cr 0&U^{-1}},
\quad U\in\calU,\ T\in\calT.$$
The group $\calU$ acts on $\calT$. For $U\in\calU$ and $T\in\calT$ we have
$T[U]:=\transpose UTU\in\calT$.
The group $\calU$ also acts on the dual lattice $\calT^*$ that consists of all
symmetric rational matrices $H$ such that $\sigma(TH)\in\gz$ for all $T\in\calT$.
This action is given by $T[\transpose U]$.
Here $\sigma$ denotes the trace.
\smallskip
Let $Y$ be a symmetric positive definite matrix. Its minimum $m(Y)$ is defined
to be the minimum value of all $Y[g]$ where $g$ runs through all
non zero integral columns. This is a continuous function on the
space of all positive definite symmetric matrices. For $C\ge 0$ we denote by
$\hz_2(C)$ the set of all symmetric complex matrices $Z$ with positive
definite imaginary part $Y$ such that $m(Y)>C$. This is an open subset
of the set of all symmetric matrices. The case
$C=0$ is the Siegel upper half-plane $\hz_2$.
The group  $\calP$ acts on
$\hz_2(C)$ through
$$\pmatrix{E&T\cr 0&E}\pmatrix{\transpose U&0\cr 0&U^{-1}}
(Z)=Z[U]+T.$$
The quotient
$$U_C=\hz_2(C)/\calP$$
is a normal complex space.
\smallskip
We are interested in the group
$$\Cl(\calP):=\lim_{C\to\infty}\Cl(U_C).$$
This group may be very big and we are only interested in a small
part of it, the Heegner part:
\smallskip
Let $S$ be a fixed matrix in $\calT^*$ with negative
determinant. For any real number $d$,
$$\{Z\in\hz_2;\quad\sigma(ZS)=d\}$$
is a non-empty set of codimension 1.
We denote by $H(S,d)$ the set of all $Z$ with
$$\sigma(ZS[U])\equiv d\; \mod 1\quad\hbox{for some}\ U\in\calU.$$
This can be considered as a $\calP$-invariant divisor where
the multiplicities are taken to be 1.
The matrix $S$ is called primitive if it can not be written in the form
$S=tS_1$ where $t>1$ is a natural number different form 1 and $S_1\in \calT^*$.
Then the equation $\sigma(ST)=1$ has a solution $T\in\calT$. Hence,
in the primitive case, $H(S,d)$ is the $\calP$-orbit of the divisor
$\sigma(ZS)=d$.
In the general case, it is a finite union
of such $\calP$-orbits. This follows from the trivial formula
$$H(NS,d)=\bigcup_{\nu=1}^NH\Bigl(S,{d+\nu-1\over N}\Bigr)\qquad
(S\ \hbox{primitive},\ N\in\nz).$$
We denote the subgroup of
$\Cl(U_C)$ spanned by these divisors by
$\Cl_{\hbox{\sevenrm Heeg}}(U_C)$. This part can be described by local
Borcherds products. Local Borcherds products have been treated in the literature
in  different contexts. In [BF], the local divisor class group of a generic
point of a one-dimensional boundary component has been treated, even
more general in the context of the group $\O(2,n)$.
Another case that has been treated by Bruinier are the cusps of
Hilbert modular surfaces in [BZ]. Here we have to consider the 0-dimensional cusps
of Siegel threefolds.
\smallskip
The convergence of the local Borcherds products simply will
rest on the following result, in which we use the notation
$e(Z)=\exp2\pii \sigma(Z)$.
\proclaim
{Lemma}
{Let $n>0$ be a natural number.
The series
$$\sum_He(-HY),\quad Y\ \hbox{symmetric positive definite},$$
converges. Here $H$ runs through all integral matrices
with determinant $\det H=-n$ and such that $\sigma(H)\ge 0$.}
KonvB%
\finishproclaim
The proof is left to the reader.\qed
\smallskip
We want to construct a holomorphic function with
divisor $H(S,d)$.
This can be done through the local Borcherds product
$$B(Z)=\prod_H\bigl(e^{2\pii\epsilon(H)(\sigma(HZ)-d)}-1\bigr).$$
Here $H$ runs through all matrices of the form
$H=S[\transpose U]$ with $U\in\calU$ and
$$\epsilon(H)=\cases{+1& if $\sigma(H)\ge 0$,\cr
-1& else.\cr}$$
From Lemma \KonvB\ follows that this product converges in $\hz_2$ and defines an analytic
function there.
The function $B$ is periodic with respect to $\calT$.
\proclaim
{Lemma}
{Let $B(Z)$ be the  local Borcherds product with zero
set $H(S,d)$.
Assume that $H(S,d)$ is in the ramification of the map
$\hz_2(C)/\calT\to\hz_2(C)/\calP$. Then the multiplicity of the zeros of
$B(Z)$ is two. This is also the ramification order.
}
MultBorch%
\finishproclaim
{\it Proof.\/} Assume that the set $\sigma(SZ)=d$ is in the ramification.
This means that there exists a substitution $Z\mapsto Z[U]+H$ in $\calP$ that
fixes it.
This transformation must be of finite order, and hence $U$ is of finite order.
The set
$$\{Y>0;\quad\sigma(SY)=0\}$$
is not empty and two-dimensional. The equation $Y[U]=Y$ holds on this set.
We claim that $\det U=-1$. This follows from the fact that an element
from $\SL(2,\gz)$ has at most one fixed point with respect to the standard action
of $\SL(2,\gz)$ on the upper half plane. The same argument shows $U^2=\pm E$.
The minus sign is not possible, since otherwise the two eigenvalues of $U$ would be
equal (both $\pm\imag$) and then $U$ would be a multiple of the unit matrix.
The known list of the elements of finite order in $\GL(2,\gz)$ shows that each
of the matrices $U$ is conjugate to one of the following two:
$$U_1=\pmatrix{1&0\cr0&-1},\quad U_2=\pmatrix{0&1\cr1&0}.$$
In the first case the solutions of $\sigma(Y[U_1])=Y$ are the diagonal matrices.
The corresponding matrix $S$ must be a constant multiple of $U_2$.
Analogously, the matrix $S$ that corresponds to $U_2$ is a constant multiple of
$U_1$. In both cases we have $S[\transpose U]=-S$. Hence this holds in general.
We also see that the ramification order is two.
\smallskip
We know
$$\sigma(SZ)=d\Longrightarrow \sigma(S(Z[U]+H)=d.$$
Hence we get
$$\sigma(SZ)=d\Longrightarrow \sigma(SZ)=\sigma(SH)-d.$$
This gives $2d=\sigma(SH)\in\gz$. Since in the Borcherds product besides $H=S$ also
$H=-S=S[\transpose U]$ occurs, we obtain that the set $\sigma(SZ)=d$ is a
at least double zero. It is easy to see that the multiplicity is really 2.
\qed
\smallskip
For $U\in\calU$ we have that
$$J(U,Z):=B(Z[U])/B(Z)$$
is a cocycle,
$$J(UV,Z)=J(V,Z[U])J(U,Z).$$
(As in group cohomology usual, we consider the action
of $\calU$ from the left. This is given by
$(U,f)\loma g$ where $g(Z)=Z[U]$.)
So $J(U,Z)$ represents an element of
$$H^1(\calU,\calO^*(\hz_2(C)/\calT)).$$
We denote the part of this group that is generated by the
$J(U,Z)$ above by
$$H^1_{\hbox{\sevenrm Heeg}}(\calU,\calO^*(\hz_2(C)/\calT)).$$
Let $B(Z)$ be a local Borcherds product.
Its zero divisor, considered on $\hz_2(C)$, is invariant under $\calP$.
Hence it induces a divisor on $U_C$.
By definition, the multiplicity of a component is the multiplicity
of $B$ considered on $\hz_2(C)$ divided by the ramification order (which is 1 or 2).
By Lemma \MultBorch\ this quotient is an integer.
Hence we obtain the following result.
\proclaim
{Lemma}
{There is a homomorphism
$$H^1_{\hbox{\sevenrm Heeg}}(\calU,\calO^*(\hz_2(C)/\calT))\lo
\Cl(U_C)$$
that attaches to the class of the cocycle $B(Z[U])/B(Z)$ the induced divisor
on $U_C$ with multiplicities as described above.}
HomBorch%
\finishproclaim
An automorphy factor coming from a local Borcherds product
is of a very simple form because of the following lemma.
\proclaim
{Lemma}
{Let $U\in\GL(2,\gz)$ and $n>0$ a natural number.
There exist only finitely many symmetric integral $T$ with
the properties
$$\det T=-n,\quad \epsilon(T)\ne\epsilon(T[U]).$$}
NurE%
\finishproclaim
The proof is left to the reader.\qed
\smallskip
The lemma implies that the automorphy factor
$J(U,Z)=B(Z[U])/B(Z)$
is of a very simple form.
We have to consider quotients of the type
$${e^x-1\over e^{-x}-1}=-e^x.$$
Hence we see: The cocycle
related to a local Borcherds product
is of the form
$$J(U,Z)=C_U e(ZH_U),$$
where $C_U$ is a constant of absolut value 1 and $H_U\in\calT^*$.
They are given as follows:
\proclaim
{Lemma}
{The cocycle of the local Borcherds product
$$B(Z)=\prod_H\bigl(e^{2\pii\epsilon(H)(\sigma(HZ)-d)}-1\bigr)$$
associated to $H(S,d)$
is of the form
$$J(U,Z):=B(Z[U])/B(Z)=C_U e^{2\pii\sigma(ZH_U)}.$$
If $\calH$ denotes the (finite) set
of all $H=S[\transpose V]$, $V\in\calU$, such that
$\epsilon(H)\ne\epsilon(H[\transpose U^{-1}])$,
then
$$H_U=-\sum_{H\in\calH} \epsilon(H)H,$$
and
$$C_U=(-1)^{\#\calH}\exp\Bigl(2\pii d\sum_{H\in\calH}\;\epsilon(H)\Bigr).$$
}
CocB%
\finishproclaim
We call the $\calU$-module
of the functions
$$C\cdot e(ZH);\quad\betr C=1,\quad H\in\calT^*,$$
by $\calE$. Then our cocycles are in the image of the natural map
$$H^1(\calU,\calE)\lo
H^1(\calU,\calO^*(\hz_2(C)/\calT)).$$
{\bf Remark.\/} We also can consider the submodule $\calE_0\subset\calE$
of all elements with $C=1$. Then $\calE_0\otimes_\gz\rz$
is isomorphic to $\Sym^2(\rz^2)$. Hence there is a link to the
Eichler cohomology group $H^1(\calU,\Sym^2(\rz^2))$ and from there
are relations to elliptic cusp forms of weight 4.
\smallskip
We denote the part of $H^1(\calU,\calE)$ coming from local Borcherds products
by $H^1_{\hbox{\sevenrm Heeg}}(\calU,\calE)$. We have proved the following result.
\proclaim
{Lemma}
{There is a natural homomorphism
$$H^1_{\hbox{\sevenrm Heeg}}(\calU,\calE)\lo\Cl(U_C).$$}
NatHom%
\finishproclaim
\neupara{A very particular case}%
We take for $\calP$ the parabolic group described in Sect.~3.
\proclaim
{Lemma}
{The group $\calU$ can be generated by the matrices
$$V_0=\pmatrix{1&0\cr 0&-1},\quad V_1=\pmatrix{1&2\cr0&1},
\quad V_2=\pmatrix{1&0\cr4&1},\quad V_3=\pmatrix{3&2\cr4&3}.$$
}
UgeN%
\finishproclaim
\vskip-\medskipamount\noindent
We skip the proof and simply mention that the group $\calU\cap\SL(2,\gz)$ is conjugate
to the group $\Gamma_0[8]$. On can use the program MAGMA [BMP] to get generators.\qed
\proclaim
{Theorem}
{In the special case of $\calP$ above, we have
$$\lim_{C\to\infty}\Cl(U_C)\cong\gz.$$}
SpezCl%
\finishproclaim
\vskip-\bigskipamount\noindent
{\it Proof.\/}
Recall (Sect.~2) that there exists a subgroup $\Gamma''\subset\Sp(2,\gz)$
of finite
index such that $\calP$ is just the
parabolic subgroup defined by $C=0$.
We use some details of the Satake compactification
$X_\Gamma''=\hz_2^*/\Gamma''$
as described in [Fr] for example.
For $C$ large enough the natural map
$$\hz_2(C)/\calP\lo\hz_2/\Gamma''$$
is an open embedding. Hence we can consider then
$U_C$ as a subset of $\hz_2/\Gamma''$.
There exists a fundamental system
of open neighborhoods $\hat U_C$ of the standard zero dimensional
boundary component $\infty$ in $X_{\Gamma''}$ such that
$$\hat U_C\cap (\hz_2/\Gamma)=U_C\qquad (C>>0).$$
The complement of $U_C$ in $\hat U_C$ is a curve.
From [FS2] we know that the germ $(X_\Gamma'',\infty)$ is isomorphic
to the germ of a quadric at the origin. In [FS2] we have taken
the quadric $Q$ defined by $z_1z_4=z_2z_3$.
Now we go back to $\hat U(C)$. The groups
$\Cl(U_C)$ and $\Cl(\hat U_C-\{\infty\})$
agree, since divisors can be extended over codimension $\ge 2$.
Now we can apply Lemma \ClNode\ to
complete the proof.
\qed\smallskip
We are interested in some very simple divisors, namely those which correspond to
the divisors $D_i^\pm$ in Sect.~2. Using the coordinates described in
Lemma \NodeTrans, we have the following result.
\proclaim
{Lemma}
{The equations of the transformed divisors $D_i^\pm$
in the model $\overline{\hz_2/\Gamma''}$ are (set theoretically):
$$\eqalign{D_1^\pm:\quad&\vartheta\Bigl[{11\atop 00}\Bigr](2Z)=
\vartheta\Bigl[{10\atop 00}\Bigr](Z)\pm\vartheta\Bigl[{10\atop 01}\Bigr](Z)=0,\cr
D_2^\pm:\quad&\vartheta\Bigl[{01\atop 00}\Bigr](Z)^2=
\vartheta\Bigl[{00\atop 00}\Bigr](Z)\vartheta\Bigl[{10\atop 00}\Bigr](Z)\pm
\vartheta\Bigl[{00\atop 01}\Bigr](Z)\vartheta\Bigl[{10\atop 01}\Bigr](Z)=0,\cr
D_3^\pm:\quad&\vartheta\Bigl[{10\atop 00}\Bigr](2Z)-\vartheta\Bigl[{01\atop 00}\Bigr](2Z)=\cr&
\vartheta\Bigl[{10\atop 00}\Bigr](Z)\vartheta\Bigl[{00\atop 10}\Bigr](Z/2)
\pm\vartheta\Bigl[{00\atop 01}\Bigr](Z)
\vartheta\Bigl[{10\atop 01}\Bigr](Z)=0.
\cr
}$$
}
TraDiv%
\finishproclaim
The proof is a straightforward calculation. We only mention that the theta
relations
$$\vartheta\Bigl[{a\atop b}\Bigr](Z)^2=
\sum_x(-1)^{\transpose xb}\vartheta\Bigl[{a+x\atop 0}\Bigr](2Z)
\vartheta\Bigl[{x\atop 0}\Bigr](2Z).$$
$$\vartheta\Bigl[{00\atop 10}\Bigr](Z/2)=\vartheta\Bigl[{00\atop 00}\Bigr](2Z)-
\vartheta\Bigl[{10\atop 00}\Bigr](2Z)+\vartheta\Bigl[{01\atop 00}\Bigr](2Z)-
\vartheta\Bigl[{11\atop 00}\Bigr](2Z)$$
have been used.\qed
\smallskip
We also can give the equations in the Siegel upper half plane.
\proclaim
{Lemma}
{Close to the standard boundary point, the divisors $D_1^+,\dots,D_3^-$ can be defined
as images of the following sets.
\smallni
\halign{\qquad$#${\rm :}\quad\hfil&$#$\qquad\hfil&$#${\rm :}\quad\hfil&$#$\quad\hfil\cr
D_1^+&2z_1+4z_2=1,&D_1^-&2z_1=1,\cr
D_2^+&z_0+3z_1+2z_2=1,&D_2^-& z_0+z_1=1,\cr
D_3^+&z_2=z_0+2,&D_3^-& z_2=z_0.\cr
}
}
EquDiv%
\finishproclaim
{\it Proof.\/} We start with the divisor $D_1^-$:
We restrict the series
$$\eqalign{
&\vartheta\Bigl[{a_1a_2\atop b_1b_2}\Bigr](2Z)=\cr&\sum
(-1)^{b_1g_1+b_2g_2}e^{2\pii(z_0(g_1+a_1/2)^2+2z_1(g_1+a_1/2)(g_2+a_2/2)+
z_2(g_2+a_2/2)^2)}\cr}$$
to $2z_1=1$. Since $2g_1g_2$ is even we get
$$e^{\pii a_1a_2/2}\;\vartheta\Bigl[{a_1\atop b_1+a_2}\Bigr](z_0)
\cdot \vartheta\Bigl[{a_2\atop b_2+a_1}\Bigr](z_2).$$
Since $\vartheta[\gotm](z)$ vanishes for the odd characteristic $\left(1\atop 1\right)$,
we get that $\vartheta\left[{11\atop 00}\right](2Z)$ vanishes along $2z_1=1$.
\smallskip
Next we restrict $\vartheta\left[{10\atop 00}\right](Z)$ to $2z_1=1$.
The result is
$$\sum e^{\pii\left(z_0(g_1+1/2)^2+(g_1+1/2)g_2+z_2g_2^2\right)}.$$
we replace $g_1$ by $-g_1-1$. Then we use that
$$(-g_1-1/2)g_2\equiv (g_1+1/2)g_2+g_2\;\mod\;2.$$
This shows that $\vartheta\left[{10\atop 00}\right](Z)$
and $\vartheta\left[{10\atop 01}\right](Z)$ agree on $2z_1=1$.
This settles $D_1^-$.
\smallskip
Next we treat $D_1^+$: We use the unimodular transformation
$Z\mapsto Z[{10\atop21}]$. This transforms the set $2z_1+4z_2=1$ to the set
$2z_1=1$. The form $\vartheta\left[{11\atop 00}\right](2Z)$ is invariant under
this transformation. Hence it vanishes also on $2z_1+4z_2=1$.
Similarly $\vartheta\left[{10\atop 00}\right](Z)$ is invariant, but
$\vartheta\left[{10\atop 01}\right](Z)$ changes its sign.
Hence $D_1^+$ has been reduced to $D_1^-$.
\smallskip
Next we treat $D_2^-$. We use the transformation $Z\mapsto Z\left[{11\atop 01}\right]$.
It transforms the set $z_0+z_1=1$ to the set $z_1=1$ and
$\vartheta\Bigl[{01\atop 00}\Bigr](Z)^2$ is transformed to
$\vartheta\Bigl[{11\atop 00}\Bigr](Z)^2$. We have already seen that this vanishes along
$z_1=1$. The series
$$\vartheta\Bigl[{00\atop 00}\Bigr](Z),\quad\vartheta\Bigl[{10\atop 00}\Bigr](Z),\quad
\vartheta\Bigl[{00\atop 01}\Bigr](Z),\quad\vartheta\Bigl[{10\atop 01}\Bigr](Z)$$
are invariant under the transformation. For the restriction to $z_1=1$ we use
that $\vartheta\left[{a_1a_2\atop b_1b_2}\right](Z)$ goes (see above) to
$$e^{\pii a_1a_2/2}\;\vartheta\Bigl[{a_1\atop b_1+a_2}\Bigr](z_0)
\cdot \vartheta\Bigl[{a_2\atop b_2+a_1}\Bigr](z_2).$$
Hence
$$\vartheta\Bigl[{00\atop 00}\Bigr](Z)\vartheta\Bigl[{10\atop 00}\Bigr](Z)-
\vartheta\Bigl[{00\atop 01}\Bigr](Z)\vartheta\Bigl[{10\atop 01}\Bigr](Z)$$
goes to
$$\vartheta\Bigl[{0\atop 0}\Bigr](z_0)\vartheta\Bigl[{0\atop 0}\Bigr](z_2)
\vartheta\Bigl[{1\atop 0}\Bigr](z_0)\vartheta\Bigl[{0\atop 1}\Bigr](z_2)-
\vartheta\Bigl[{0\atop 0}\Bigr](z_0)\vartheta\Bigl[{0\atop 1}\Bigr](z_2)
\vartheta\Bigl[{1\atop 0}\Bigr](z_0)\vartheta\Bigl[{0\atop 2}\Bigr](z_2)$$
which is zero. This settles $D_2^-$.
\smallskip
To treat $D_2^+$ we use the transformation
$Z\mapsto Z\left[{1\phantom{-}0\,\atop 2\,-1}\right]$.
It transforms the set $z_0+z_1=1$
to the set $z_0+3z_1+2z_2=1$. It leaves
$\vartheta\left[{01\atop 00}\right](Z)^2$ invariant. Hence this function vanishes
along $z_0+3z_1+2z_1=1$. We have seen already that it changes the sign of
$\vartheta\left[{10\atop 01}\right](Z)$. This shows the transformation of the
expression
$\vartheta\left[{00\atop 00}\right](Z)\vartheta\left[{10\atop 00}\right](Z)-
\vartheta\left[{00\atop 01}\right](Z)\vartheta\left[{10\atop 01}\right](Z)$
just causes a change of the minus sign to a plus sign.
Now $D_2^\pm$ is settled.
\smallskip
We treat $D_3^-$.
First one sees that the transformation $Z\mapsto Z\left[{01\atop10}\right]$
interchanges
$\vartheta\bigl[{10\atop 00}\bigr](2Z)$ and $\vartheta\bigl[{01\atop 00}\bigr](2Z)$.
Since this transformation is the identity on $z_0=z_2$, the difference vanishes
on this set.
\smallskip
We still have to investigate the four series
\smallni
\halign{\qquad$#$\hfil&\quad$#$\hfil&\qquad$#$\hfil&\quad$#$\hfil\cr
A:&\vartheta\bigl[{10\atop00}](Z),&B:&\vartheta\bigl[{00\atop10}](Z/2),\cr
C:&\vartheta\bigl[{00\atop01}](Z),&D:&\vartheta\bigl[{10\atop01}](Z).\cr}
\smallskip
The claim is that $AB=CD$ on $z_0=z_2$.
We use the transformation $Z\mapsto Z\bigl[{1\phantom{-}1\atop1-1}\bigr]$.
Then the set $z_0=z_2$ transforms to $z_1=0$.
We will use the formula
$$Z\left[\Bigl({1\phantom{-}1\atop1-1}\Bigr)\Bigl({g_1\atop g_2}\Bigr)\right]
=z_0(g_1+g_2)^2+z_2(g_1-g_2)^2\quad\hbox{for}\quad z_1=0.$$
If $g_1,g_2$ runs through all integers than $h_1=g_1+g_2$, $h_2=g_1-g_2$
runs through all pairs of integers such that $h_1+h_2$ is even.
For $A$ we get the expression
$$\sum_{h_1+h_2\ \hbox{\sevenrm even}}
e^{\pii(z_0(h_1+1/2)^2+z_2(h_2+1/2)^2)}.$$
We divide this sum into two parts where $h_1,h_2$ both are even or both are odd.
Both partial sum are equal and we get:
$$A:\quad 2\vartheta\Bigl[{1/2\atop0}\Bigr](4z_0)
\vartheta\Bigl[{1/2\atop 0}\Bigr](4z_2).$$
In a similar way we we split $B,C,D$. The result is
$$\eqalign{
B:\quad  &\vartheta\Bigl[{0\atop 1}\Bigr](2z_0)
\vartheta\Bigl[{0\atop 1}\Bigr](2z_2),\cr
C:\quad  &\vartheta\Bigl[{0\atop 1}\Bigr](4z_0)
\vartheta\Bigl[{0\atop 1}\Bigr](4z_2),\cr
D:\quad  &2e^{-\pii/4}\vartheta\Bigl[{1/2\atop 1}\Bigr](4z_0)
\vartheta\Bigl[{1/2\atop 1}\Bigr](4z_2).\cr}$$
The relation $AB=CD$ means:
$$\vartheta\Bigl[{1/2\atop 0}\Bigr](2z) \vartheta\Bigl[{0\atop 1}\Bigr](z)
=e^{-\pii/4}\vartheta\Bigl[{0\atop 1}\Bigr](2z)\vartheta\Bigl[{1/2\atop 1}\Bigr](2z).$$
This relation between theta series of one variable is left as an exercise.
\smallskip
The case $D_3^+$ is similar.\qed
\proclaim
{Lemma}
{The divisors  $D_i^-$ correspond (close to the standard node)
to divisors of the form
$H(S,d)$:
$$\eqalign{
&D_1^-:\qquad S=1/4\pmatrix{0&1\cr 1&0},\quad d={1\over 4},\cr
&D_2^-:\qquad S=1/4\pmatrix{2&1\cr 1&0},\quad d={1\over 2},\cr
&D_3^-:\qquad S=1/4\pmatrix{1&0\cr 0&-1},\quad d=0.
\cr}$$}
AlsHeeg%
\finishproclaim
Notice that $S$  is a primitive element of $\calT^*$.
We compute the value $J(U,Z)$ of the cocycle of the associated
local Borcherds product for some $U$. Recall (\CocB) that
$$J(U,Z):=B(Z[U])/B(Z)=C_U e^{2\pii\sigma(ZH_U)}.$$
The constant $C_U$ is not important, since it is a root of unity that disappears if
one takes a suitable power. Hence we only give our attention to $H_U$.
We compute them for the generators given in Lemma \UgeN.
Evaluation of
the formula given in Lemma \CocB\ gives the following values for the $8H_{V_i}$:
$$\eqalign{&
D_1^-:\qquad\pmatrix{0&0\cr0&0},\quad\pmatrix{-4&-2\cr-2&0},\quad
\pmatrix{0&-2\cr-2&-8},\quad\pmatrix{-12&-18\cr-18&-24}.\cr&
D_2^-:\qquad\pmatrix{0&0\cr0&0},\quad\pmatrix{-4&-2\cr-2&0},\quad
\pmatrix{0&0\cr0&0},\quad\pmatrix{0&0\cr0&0}.\cr&
D_3^-:\qquad\pmatrix{0&0\cr0&0},\quad\pmatrix{-4&-2\cr-2&0},\quad
\pmatrix{0&0\cr0&0},\quad\pmatrix{0&0\cr0&0}.\cr}$$
To compare the we make use of trivial divisors
$e^{\sigma(UZ\transpose U)}/e^{\sigma (Z)}$.
For $U=V_1$ and $U=V_2$ one gets the following values
$$\eqalign{&
(V_1):\qquad\pmatrix{0&0\cr0&0},\quad\pmatrix{0&0\cr0&0},\quad
\pmatrix{0&4\cr4&16},\quad\pmatrix{8&12\cr12&16}.\cr&
(V_2):\qquad\pmatrix{0&0\cr0&0},\quad\pmatrix{4&2\cr2&0},\quad
\pmatrix{0&0\cr0&0},\quad\pmatrix{4&6\cr6&8}.\cr}$$
From these data follows the following proposition.
\proclaim
{Proposition}
{Consider the divisors $D_1^\pm,D_2^\pm,D_3^\pm$. Their classes
in $\Cl(\calX,\eta)$ are non zero and
satisfy the relations $D_i^+=-D_i^-$ and
$$D_1^\pm=-D_2^\pm=-D_3^\pm.$$}
RelCl%
\finishproclaim
\vskip-\bigskipamount\noindent
The stabilizer of a node $a$ acts on the local Picard group $\Cl(\calX,a)$.
The subgroup of index two that fixes the rulings acts as identity.
Hence we see the following result.
\proclaim
{Remark}
{By means of Proposition \RelCl\ and Theorem \RulStan\ it is possible to compute
the image of each of the $188$ basic divisors in $\Cl(\calX,\eta)$.
Since they generate $\Cl(\calX)\otimes_\gz\qz$ we have a complete description of the map
$$\Cl(\calX)\otimes_\gz\qz\lo \Cl(\calX,\eta)\otimes_\gz\qz\qquad(\cong\qz).$$}
CompComp%
\finishproclaim
Hyperplane sections of $\calX$ define line bundles and hence
get trivial in the local divisor class
groups at nodes. Hence the one dimensional representation is in the kernel of the
map
$$\Cl(\calX)\otimes_\gz\qz\lo\bigoplus_{a\;\hbox{\sevenrm node}}
\Cl(\calX,a)\otimes_\gz\qz$$
From the above explicit description one can compute that the one dimensional representation
is the precise kernel. This might be a general phenomenon
for Siegel threefolds (compare [BF], Theorem 5.4).
\neupara{Projective Resolutions}%
If $D$ is an effective divisor on an irreducible normal complex space
$X$, then we associate the ideal sheaf $\calI(D)$ of all holomorphic
functions on open subsets that satisfy $(f)\ge D$ on this subset.
The blow up of $D$, by definition, is the blow up of the ideal
sheaf $\calI(D)$.
The
blow up of an effective divisor depends only on its image in
the divisor class group $\Cl(X)$.
\smallskip
The relation $\calI(D)^n=\calI(nD)$ is not true
in general. But in our situation it is true.
\proclaim
{Lemma}
{Let $Q$ be the affine quadric defined by $z_1z_4=z_2z_3$
and let $P$ be the prim divisor defined by $z_1=z_2=0$.
Then
$$\calI(nP)=\calI(P)^n.$$}
nPPn%
\finishproclaim
{\it Proof.\/} We have to show $\calI(nP)\subset\calI(P)^n$ for $n>0$
since
the converse inclusion is trivial. Hence we have to consider
an element $f$ of the (analytic) local ring $\calO_{Q,0}$ that vanishes
along $P$ of order at least $n$. We can take this element as
the restriction of a power series
$$F\in R:=\cz\{z_1,z_2,z_3,z_4\}.$$
We expand $F$ as power series of the variables $z_1,z_2$,
$$F=\sum g_{i,j}(z_3,z_4)z_1^i z_2^j.$$
Using the relation $z_1z_4=z_2z_3$ we can modify $F$ in such a way that
$g_{i,j}$ is independent of $z_3$ in the case $j>0$.
Now we claim that $g_{i,j}=0$ if $i+j<n$.
The proof is given by induction on $k=i+j$. The beginning of the induction
($k=0$ is trivial. To get the idea, we treat the case
$k=1$ (assuming $n\ge 2$).
We know that the restriction of
$$g_{1,0}z_1+g_{0,2}z_2$$
to $Q$ vanishes along $z_1=z_2=0$ in at least second order.
We test this in the chart $z_3\ne 0$. The divisor is given
in this chart by one equation $z_1=0$. The claim is that
$$g_{1,0}(z_3,z_4)z_1+g_{0,1}(z_4)z_1z_4/z_3$$
vanishes of at least second order on $z_1=0$.
Then
$$g_{1,0}(z_3,z_4)+g_{0,1}(z_4)z_4/z_3$$
still must vanish along $z_1=0$. Since it is independent of
$z_1$ it must vanish. This gives
$$g_{0,1}(z_4)z_4=-g_{1,0}(z_3,z_4)z_3.$$
Since the left-hand side is independent of $z_3$
(due to our normal form) we get $g_{0,1}=0$ and then
$g_{1,0}=0$.\qed
\smallskip
Let  $(X,a)$ be a three dimensional nodal singularity.
The blow up of a divisor $D$ gives a small resolution of
the node if and only if
$D$ runs non-trivially into the node in the following sense:
the image of $D$ in the local divisor class group $\Cl(X,a)$ is not
zero. This follows easily from Lemma \nPPn\ in connection with the
the structure theorem $\Cl(X,a)\cong\gz$.
\smallskip
Let now $G\subset\bar\calH$ be s subgroup. We recall that we defined
at the end of the introduction two sets $\calA$, $\calB$ of nodes.
We reformulate the theorem at the end of the introduction in the following way.
\proclaim
{Proposition}
{Let $G\subset\bar\calH$ be a subgroup
such that $\calA\cup\calB$ is the set of all
nodes. Assume that for each node
$a\in\calB$, $a\not\in\calA$, the map
$$(\Cl(\calX)\otimes_\gz\qz)^G\lo\Cl(\calX,a)\otimes_\gz\qz$$
is not the zero map. Then there exists a projective resolution $\tilde\calX\to\calX$
such that the action of $G$ extends to a group of biholomorphic transformations
of $\tilde\calX$ and such that $G$-$\Hilb(\tilde\calX)$ is a resolution
of $\calX/G$ in the form of a projective Calabi--Yau manifold.}
HauptS%
\finishproclaim
{\it Proof.\/} First we blow up the nodes from $\calA$ to get a partial resolution
$\calX_1\to\calX$. The group $G$ extends since the set $\calA$ is $G$-invariant.
If $\calX_1$ is smooth we are done. Otherwise
we choose some node $a\in\calX$ that remains singular in $\calX_1$.
By assumption
there exists a divisor $D$ such that its class in $\Cl(\calX)\otimes_\gz\qz$
is $G$-invariant and that
is non-trivial at $a$. We can assume that $D$ is effective, since we can
add to $D$ the divisor of a homogenous polynomial in $Y_0,\dots,X_3$.
We can also assume that $D$ itself (and not only its class)
is $G$-invariant, since we can replace $G$ by the sum of all $g(D)$, $g\in G$.
Now we want to blow up $D$.
To be precise, we take the transform $D_1$ of $D$ in $\calX_1$,
and blow up $D_1$ to produce
a partial resolution $\calX_2\to\calX_1$.
Of course the blow up doesn't change anything in the smooth
part of $\calX_1$ because there any divisor is locally principal.
The resolution $\calX_2\to\calX_1$ is a small resolution of certain nodes
including the node $a$. If $\calX_2$ is already smooth we are done.
Otherwise we proceed in the same manner.
\smallskip
The constructed $\tilde\calX$ has the property that for each point of the quotient
$\tilde\calX/G$ there exists a small (analytic) neighborhood $U$ and a
holomorphic 3-form without zeros on the regular locus of $U$.
For the nodes $a\in\calA$ one has to observe the following. The Calabi--Yau form
gets a zero of order one along the exceptional divisor of the blow up.
But, by the definition of $\calA$ there is an (central) element $g\in G_a$ that fixes
the exceptional divisor pointwise.
Due to the ramification, the Calabi--Yau form, considered on the quotient
of the blow up by $g$, doesn't vanish along the image of the exceptional divisor.
\smallskip
As a consequence, all singularities of $\tilde\calX/G$
are of the form $\cz^3/H$, where $H\subset\SL(3,\cz)$, is a finite group.
But then the results of [BKS] show that the $G$-Hilbert scheme
$G$-$\Hilb(\tilde\calX)$ gives a projective resolution in the form of a
Calabi--Yau manifold
(compare [FS1],  Theorem 2.6).
\qed
\neupara{Examples}%
1) It can happen that a group $G\subset\bar\calH$ contains the conjugates of
the transformation
$(Y_0,-Y_1,Y_2,-Y_3,X_0,X_1,-X_2,-X_3)$.
For such a group
all nodes are contained in $\calA$. Hence
$\calX/G$ admits a projective Calabi-Yau resolution.
These examples has been described in detail in [FS2]
(see Lemma 8.1, Remark 8.3 and Corollary 8.4).
\medni
2) We start with a counter example. There exist groups $G\subset\calG$ of order
32 that acts freely on $\calX$. An example is the group generated
by the following transformations.
$$\eqalign{
&(-Y_1,-\imag Y_0,Y_3,-\imag Y_2,-\imag X_1,X_0,\imag X_3,X_2),\cr
&(-\imag Y_2,Y_3,Y_0,\imag Y_1,X_3,-X_2,-\imag X_1,-\imag X_0),\cr
&(-Y_0,Y_1,Y_2,-Y_3,X_0,-X_1,X_2,-X_3),\cr
&(\imag Y_0,\imag Y_1,\imag Y_2,\imag Y_3,-\imag X_0,-\imag X_1,-\imag X_2,-\imag X_3),\cr
&(Y_0,Y_1,-Y_2,-Y_3,X_0,X_1,-X_2,-X_3),\cr
&(\imag Y_0,\imag Y_1,\imag Y_2,\imag Y_3,\imag X_0,\imag X_1,\imag X_2,\imag X_3)\cr}
$$
In this case all nodes are of type $\calB$. But the space $(\Cl(\calX)\otimes_\gz\qz)^G$ is
one dimensional (generated by a hyperplane section). The image in
$\Cl(\calX,a)\otimes_\gz\qz$ is zero for all nodes $a$.
We claim even more. There is no projective Calabi--Yau manifold $M$
that is birational
equivalent to $\calX/G$. This follows from a theorem of Kollar that states that
two bimeromorphic equivalent models are related by flops
([Ko], Theorem 4.9). Hence $M$ must be obtained
from $\calX/G$ by a small resolution. The pull back to $\calX$ would give a
projective $G$-invariant resolution.
It is easy to show that this is the blow up of a $G$-invariant divisor.%
\footnote{*)}{\ninepoint
We are very grateful to van Geemen who pointed out to us that there exists no
projective model for this group $G$.}
\smallskip
Using calculator we got the following result.
\proclaim
{Theorem}
{There are $54$ conjugacy classes of subgroups $G\subset\bar\calH$ that act freely
on $\calX$. Their orders are in $\{1,2,4,8,16,32\}$.
For those of order $32$ there is no projective Calabi--Yau model for
$\calX/G$ (but a weak Calabi-Yau model). Those of order $1$, $2$ and $4$ all
admit projective Calabi--Yau models.
In the case of order 8, 13 cases have a
projective  C-Y model and 7 cases do not have it.
There is one class of order $16$ with a projective Calabi--Yau model, the
other $12$ classes do not have it.}
freieG%
\finishproclaim
The quotient of a rigid projective manifold $M$ by a finite freely acting group is rigid.
The Euler number is $e(M/G)=e(M)/\#G$. Hence we get rigid Calabi--Yau manifolds
with Euler numbers $4,8,16,32$.
\smallskip
The
freely acting group of order $16$ with a projective Calabi--Yau model is of special
interest.
\proclaim
{Theorem}
{The group $G$
generated by the two transformations
$$\eqalign{&(-Y_2,\imag Y_3,\imag Y_0,Y_1,-\imag X_3,-\imag X_2,X_1,-X_0),\cr&
(\imag Y_1,Y_0,-\imag Y_3,Y_2,X_1,\imag X_0,X_3,-\imag X_2)}$$
has order $16$ and  acts freely on $\calX$. The quotient $\calX/G$
has a resolution in the form of
a rigid Calabi--Yau manifold  ($h^{12}=0$) with Euler number $e=4$
and Picard number $h^{11}=2$.}
freiCM%
\finishproclaim
3) We consider the group $G$ of order 2 that is generated by the involution
$$\sigma_1:\quad (Y_0,Y_1,Y_2,Y_3,-X_0,-X_1,-X_1,-X_3).$$
It generates a normal subgroup $G$ of $\bar\calG$.
The involution is fixed point free. Hence all nodes are of type
$\calB$ (and there are no ones of type $\calA$).
The divisor $D_2^+$  ($X_0X_2+X_1X_3=Y_0Y_1\pm Y_2Y_3=0$)
is invariant under $G$ and it
is non-trivial at the standard node (see Proposition \RelCl).
Since $G$ is normal we obtain by transformation for each node
a divisor in that is $G$-invariant and non-trivial at this node.
Hence we get a projective model.
\medni
4) Next we consider the group of order 2 that is generated by the involution
$$\sigma_2:=(Y_0, -Y_1, -Y_2, Y_3, X_0, -X_1, -X_2, X_3).$$
This substitution has 6 conjugates in $\bar\calH$, namely
$$\eqalign{
(Y_0, -Y_1, Y_2, -Y_3, -X_0, -X_1, X_2, X_3),\
&(Y_0, -Y_1, -Y_2, Y_3, X_0, -X_1, -X_2, X_3),\cr
(Y_0, Y_1, -Y_2, -Y_3, -X_0, X_1, -X_2, X_3),\
&(Y_0, -Y_1, -Y_2, Y_3, -X_0, X_1, X_2, -X_3),\cr
(Y_0, Y_1, -Y_2, -Y_3, X_0, -X_1, X_2, -X_3),\
&(Y_0, -Y_1, Y_2, -Y_3, X_0, X_1, -X_2, -X_3).\cr}$$
Again the divisor $D_2^+$ is invariant under them and not trivial at the standard
node. So we get projectivity again.
\medskip
There are 10 conjugacy classes of
involutions in $\bar\calH$. In [FS2] we listed them in Proposition 7.6
and computed the divisor class and Euler numbers for them.
As in the case $\sigma_1$, $\sigma_2$ above one can verify in each of the cases that
\HauptS\ applies to obtain a projective Calabi--Yau manifold again.
Then we obtain
\proclaim
{Theorem}
{The following table
describes fixed point sets of the involutions $\sigma_i$, $i\le i\le 10$,
on $\calX$ and the Hodge numbers numbers of a (projective) Calabi--Yau model
of the quotient $\calX/\sigma_i$.
\smallni\vbox{
\halign{\qquad\qquad$\sigma_{#}$:\quad\hfil&\hbox{\rm #}\qquad\hfil
&$#$\quad\hfil&$#$\hfil\cr
\omit&fixed points&h^{11}&h^{12}\cr
\noalign{\vskip1mm}
1& empty set&16&0\cr
2& 16 nodes&40&0\cr
3& 4 elliptic curves&20&4\cr
4& empty set&16&0\cr
5& empty set&16&0\cr
6& 8 conics in  planes ($\cong \pz^1$)&28&0\cr
7& 8 lines ($\cong \pz^1$) &28&0\cr
8& 2 elliptic curves&18&2\cr
9& 2 elliptic curves&18&2\cr
10&4 conics in planes ($\cong \pz^1$)&22&0\cr
}}
}
FixSet%
\finishproclaim
In the cases where elliptic curves are in the fixed point set, we get
{\it non rigid Calabi--Yau manifolds.}
\medni
5) Next we consider  a group of order 3.
\proclaim
{Proposition}
{There is only one conjugacy class of elements of order three
in $\bar\calH$. It can be represented by
$$(Y_0,Y_3,Y_1,Y_2,X_0,X_3,X_1,X_2).$$
Its fixed point set in $\calX$ doesn't contain a node.
It consists of one elliptic curve and $4$ isolated points.
The fixed point set in $\tilde\calX$ is the same.
The subgroup of order three leads to a projective
Calabi--Yau model with Hodge numbers
$$(h^{11},h^{12})=(18,2).$$
}
Fixthree%
\finishproclaim
{\it Proof (sketch).\/} We omit the computation of the fixed points.
We can extend the group to a small resolution $\tilde\calX$
of $\calX$. Now
the Euler number easily can be computed by means of the stringy formula
(compare [FS2]. The result is $e=32$.
Next we compute the Picard number $h^{11}$.
The invariant part $\dim\Cl(\calX)^G$ has rank
$12$ as can be proved  by means of the known character of the action
of $\bar\calH$ (see [FS2]). One can check that over
isolated points there is one exceptional divisor but over the elliptic
curve there are two.
Hence we have 6 exceptional divisors. We obtain $h^{11}=18$.
\qed
\neupara{A Special class of groups}%
This section rests highly on computer calculations.
We consider all subgroups of $\calH$ that are isomorphic to $(\gz/2\gz)^m$.
There are 165 conjugacy classes. In all 165 cases the group extends to a
small resolution $\calX$ such that a weak Calabi--Yau model can be
obtained by a resolution of $\calX/G$.
There are 144 classes that admit a
projective Calabi--Yau model that can be obtained as described in the introduction.
This means that the group $G$ extends to a not necessarily minimal
projective resolution
$\tilde\calX$ such the $\tilde\calX/G$ admits a resolution in the form
of a (projective) Calabi--Yau manifold. It is known that the
Euler number and divisor class number
are bimeromorphic invariant for weak Calabi--Yau manifolds (since they are related
by flops.) Therefore we always can use a small minimal resolution $\tilde\calX$
(also a non-projective one)
for the computation of the Euler and divisor class number.
This actually means that we sometimes
compute the Hodge numbers for a group with
projective Calabi--Yau model using a non-projective model.
\zwischen{Computation of the divisor class number}%
So let $G\subset\bar\calH$ be a subgroup that is isomorphic to $(\gz/2\gz)^m$.
We choose a small resolution $\tilde\calX$ such that $G$ extends.
We have to determine the divisor class number of a Calabi--Yau resolution of
$\tilde\calX/G$. This is the sum of the divisor class number of
$\calX/G$ and the number of exceptional divisors.
The divisor class number of $\calX/G$ can be computed since we know the
representation of $\bar\calH$ on $\Cl(\calX)\otimes_\gz\qz$.
This can be done by a program. To get the number of exceptional divisors
we need information about the singularities of $\calX/G$.
Locally they are of the form $\cz^3/H$ with an abelian group $H\subset\SL(3,\gz)$.
where all elements of $H$ have order $\le 2$.
After diagonalization $H$ consists of sign changes. Since the determinants
are one we have an even number of sign changes. In this cases the
resolution is easy to produce (see [FS1]). The result is that the singular locus
of $\cz^3$ consists of one line or of two crossing lines. In the first case
we have one exceptional divisor in the resolution in the second case we have two.
This shows that the number of exceptional divisors of a (weak)
Calabi-Yau resolution of $\tilde\calX/G$ equals the number of components of the
fixed point locus. We want to express this in the singular model $\calX$.
So let $C\cong \pz^1$ be an exceptional curve that projects to a node in $\calX$.
Then there is an $g\in G$, $g\ne e$, that fixes the node.
This node must be an isolated fixed point of $g$ since otherwise the fixed locus
on $\tilde\calX$ would be not smooth. So we see the following lemma.
\proclaim
{Lemma}
{Let $G\subset\bar\calH$ be a subgroup of type $(\gz/2\gz)^m$.
Let $a$ be the number of $G$-equivalence classes of irreducible curves in
$\calX$ that are in the fixed point locus and let $b$ be the number of
$G$-equivalence classes of nodes that are isolated fixed points.
Then the divisor class number of a weak Calabi--Yau model of $\calX/G$
equals
$$a+b+\dim (\Cl(\calX)\otimes_\gz\qz)^G.$$}
ClzW%
\finishproclaim
{\bf Computation of the Euler number}%
\medni
Our main tool will be the string theoretic formula [Ro]:
\smallni
{\it Let $X$ be a weak Calabi--Yau three fold and $G$ a finite group of biholomorphic
transformations that leave the Calabi--Yau three form invariant.
Assume  that the stabilizers
$G_a$ are contained in the group $\SL()$
of the tangent space. Then the Euler number of a
resolution of $X/G$ in the form of a Calabi--Yau manifold equals
$$e(M,G)={1\over\#G}\sum_{gh=hg}e(M^{\langle g,h\rangle}).$$
Here $M^{\langle g,h\rangle}$ denotes the common fixed point set
of $g,h$.
}
\smallni
Let $G\subset\bar\calH$ be a subgroup that extends to a small resolution
$\tilde\calX$. Then the string theoretic formula applies. Recall that the
fixed point locus of $G$ on $\tilde\calX$ is a curve $C$
with smooth irreducible components. Let $S$ be the singular locus of $S$
(crossing points). Then $e(C)=e(C-S)+\#S$.
The Euler number of $C-S$ is the sum of the Euler numbers of its connected components.
The Euler number of a smooth non-compact curve $K$ with compactification $\bar K$ is just
$$e(K)=2-2g(\bar K)-\#(\bar K-K)\qquad (g(\bar K)=\hbox{genus}).$$
Since we know the equations of $\calX$ it is no problem to compute
the fixed point loci on $\calX$. This is the image of $C$ in $\calX$.
We know already that the fixed point loci of elements of order two consist
of (smooth) rational and elliptic curves and isolated points (nodes).
Hence we have just to analyze what happens on a exceptional $\pz^1$ in $\tilde\calX$.
This is clarified by the following two lemmas.
\proclaim
{Lemma}
{Let $g\in\bar\calH$ an element of order two which fixes a node $a$.
There are two cases:
\vskip1mm
\item{\rm1)} $a$ is an isolated fixed point of $g$.
In this case $g$ extends as identity on the exceptional $\pz^1$ over $a$.
\item{\rm 2)} $a$ is not an isolated fixed point of $g$. Then $g$ has precisely
two fixed points on the exceptional $\pz^1$.
These two fixed pointy are intersection points of two further fixed point
curves of $g$ which are visible already in $\calX$.
\vskip0pt}
ordZw%
\finishproclaim
This  lemma follows from our investigations of the involutions
during the proof of \FixSet.\qed
\proclaim
{Lemma}
{Let $g$, $h$ be two different commuting elements of $\bar\calH$
of order two. Assume that they fix a joint node
$a$. Then there are two possibilities for
the joint fixed point locus on the exceptional $\pz^1$ over $a$.
\vskip1mm
\item{\rm1)} It consists of two points. This happens if $a$ is
an isolated fixed point of one of the three $\{g,h,gh\}$.
\item{\rm2)} It is empty. This happens if it is not an isolated
fixpoint of any of the three $\{g,h,gh\}$.
\vskip0pt
}
DifEmp%
\finishproclaim
{\it Proof.\/}
We first mention that it cannot happen that $a$ is an isolated
fixed point of two of the $g$, $h$, $gh$.
The reason is that there is only one conjugacy class of elements
of order two which fix a node as isolated fixed point
(the conjugacy class of $\sigma_2$). This conjugacy class
consists of 6 elements and each of them fixes
16 nodes and nothing else. The 6 blocks of 16 nodes are pairwise
disjoint and exhaust all $96=6\cdot 16$ nodes.
\smallskip
Case 1) is clear, assume for example that $gh$ has the node as
an isolated fixed point. Then it acts as identity
on the exceptional fibre. Then $g$ and $h$ are inverse
on this fibre and the same fixed points. From \ordZw\ we know that they
are two.
\smallskip
We treat the second case.
From the assumption follows that $g$ and $h$ are different.
We know that $g$ has precisely two fixed points on the exceptional
fibre over $a$. We choose
the biholomorphic map $\pz^1\cong\bar\cz$ such that $g$ acts as $g(z)=-z$.
Since $h$ commutes with $g$ it acts on the two fixed points.
It cannot fix both since then $h$ would we equal to $g$.
Hence $h$ permutes $0$ and $\infty$ and the only possibilities
are the transformations $h(z)=\pm 1/z$.
But they have different fixed points. Hence the joint fixed point set
on the $\pz^1$ is empty.\qed
\smallskip
We have collected all what we need for the computation of the Euler number.
A computer calculation gives now the following list.
\medni
There are 40 different pairs $(\hbox{cl},e)$ of divisor class numbers
and Euler numbers of weak Calabi--Yau manifolds
that are produced by subgroups of $\bar\calH$ isomorphic
to $(\gz/2\gz)^m$.
$$\eqalign{&(28, 56), (20, 40), (16, 32), (14, 28), (15, 26), (12, 16), (10, 20),
(8, 16), (4, 8),\cr&  (26, 50), (10, 8), (22, 44), (14, 26), (18, 28),
(13, 20), (44, 88), (11, 16), (16, 28),\cr& (34, 68), (20, 32), (16, 16),
(14, 20), (13, 26), (22, 40), (15, 28), (41, 82), (26, 52),\cr& (6, 8),
(12, 8), (17, 32), (19, 38), (70, 140), (14, 16), (12, 20), (10, 16),
(9, 14),\cr& (46, 92), (40, 80), (32, 64), (18, 32).\cr}$$
\smallni
There are 33 different pairs $(h^{11},h^{12})$ of Hodge numbers
numbers which are produced by subgroups of $\bar\calH$ isomorphic
to $(\gz/2\gz)^m$ and with a projective model as described in the introduction.
$$\eqalign{&(14,0), (12,4), (20,0), (14,1), (6,2), (26,0), (32,0), (26,1), (70,0),
(18,2),\cr& (16,8), (44,0), (15,2), (15,1), (12,2), (10,2), (14,6), (9,2),
(4,0), (22,2),\cr& (10,0), (16,0), (14,4), (12,8), (22,0), (20,4), (13,0),
(28,0), (19,0), (16,2),\cr& (34,0), (40,0), (8,0).\cr}$$
The list seems to contain several examples that are not contained
in the (physicists) literature (see [CD]).
\neupara{Two more  examples}%
We consider the group $H$ of order 16 that we described in Theorem \freiCM.
There are transformation of order 3 in $\bar\calH$ that normalize this group,
for example
$$h_3:\qquad(Y_0, -Y_3, \imag Y_1, \imag Y_2, X_3, X_1, X_0, X_2).$$
We consider now the group $H_{48}$ of that is generated by this element and the group
in Theorem \freiCM. The order of this group is 48.
Using Proposition \HauptS\ on can check that $H_{48}$ extends to a projective
small resolution $\tilde\calX\to\calX$.
All 32 elements
which are not contained in $H$ are of order 3. The 16 subgroups of order 3 of
$H_{48}$ are conjugated under $H$. This means that the fixed point locus of
$h_3$ in $\calX$ maps under $\calX\to\calX/H$ biholomorphic to the whole
fixed point locus of $H_{48}/H$. Recall the the fixed point locus of
$h_3$ does not contain a node and consists of an elliptic curve and 4 points.
The $H_{48}$-invariant part of $\Cl(\calX)$ has rank two. Hence the Picard number
of a Calabi--Yau resolution of $\calX/H_{48}$ equals $8$.
The string theoretic formula applied to $\tilde\calX/H,H_{48}/H$ gives
$e=(1/3)(4+8\cdot 4)=12$. Hence we obtain a Calabi--Yau manifold with
Hodge numbers
$$h^{11}=8,\quad h^{12}=2\qquad(e=12).\leqno H_{48}:$$
Our last example starts with the same group $H$ but we extend it now by a involution
$$h_2:\qquad (Y_2, \imag Y_3, -\imag Y_0, -Y_1, -\imag X_3, \imag X_2, -X_1, X_0)$$
This involution normalizes $H$. Hence we get a group $H_{32}$ of order $32$.
Using Proposition \HauptS\ on can check that $H_{32}$ extends to a projective
small resolution $\tilde\calX\to\calX$.
The whole fixed point locus of $H_{32}$ consists of 24 elliptic curves and 32
rational curves. The image in $\calX/H$ (i.e.~the fixed poind locus of $h_2$
acting on $\calX/H$ consists of 3 elliptic curves and 2 rational curves.
One obtains a (projective) Calabi--Yau manifold with Hodge numbers
$$h^{11}=7,\quad h^{12}=3\qquad(e=8).\leqno H_{32}:$$
{\bf Final Remark.} 
With a refinement of the described methods we could compute  $1\,344$ of the
$4\,117$ cases the Euler number where $786$ admit a projective Calabi--Yau model.
The Euler numbers of weak Calabi--Yau models that we got are
$$\eqalign{&
2, 4, 6, 8, 10, 12, 14, 16, 18, 20, 22, 24, 26, 28, 30, 32, 34, 38, 40, 44,
46, 50, 52, 56,\cr& 58, 64, 68, 70, 80, 82, 88, 92, 100, 140.\cr}$$
Euler numbers of a projective Calabi--Yau model from this list are
$$4, 8, 12, 14, 16, 20, 24, 26, 28, 32, 38, 40, 44, 50, 52, 56, 64, 68, 80, 88,
100, 140.$$
The Euler numbers, we could get, are all positive. This is clear from the
formula for string theoretic
Euler number in case that the group extends to a small resolution.
So we cannot get mirror pairs.
\vskip1cm\noindent
{\paragratit References}%
\bigni
\item{[BMP]} Bosma,\ W., Cannon,\ J., Playoust,\ C.:
{\it The Magma algebra system. I. The user language,\/}
J. Symbolic Comput. {\bf24} (3-4), 235--265 (1997)
\medskip
\item{[BKR]}
Bridgeland, T., King, A., Reid, M.:
{\it The McKay correspondence as an equivalence of derived categories,\/}
J.~Amer.~Math.~Soc. {\bf 14}, 535–-554 (2001)
\medskip
\item{[BF]} Bruinier,\ J.H., Freitag,\ E.: {\it Local Borcherds products}
Annales de L'Institute Fourier Grenoble, Tome {\bf 51}, 1-26 (2001)
\medskip
\item{[BZ]} Bruinier, J.H., van der Geer, G., Harder,G., Zagier,D.:
{\it The 1-2-3 of modular forms,\/}
Universitext. Springer, Berlin, Heidelberg (2008)
\medskip
\item{[CD]} Candelas, P., Davis, R.: {\it New Calabi--Yau manifolds with small
Hodge numbers,\/} Fortschritte der Physik, Vol.\ {\bf 58}, Issue 4--5,
383--466 (2010)
\medskip
\item{[CFS]} Cynk, S., Freitag, E., Salvati--Manni, R.:
{\it The geometry and arithmetic of a
Calabi-Yau  Siegel threefold,\/}
accepted by International Journal of Mathematics,
arXiv: 1004.2997 (2010)
\medskip
\item{[CM]} Cynk, S., Meyer, C.:
{\it Modular Calabi--Yau threefolds of level eight,\/}
Internat.~J.~Math.\ {\bf 18}, no.~3, 331--347 (2007)
\medskip
\item{[Fr]} Freitag, E.: {\it Siegelsche Modulfunktionen,} Grundlehren
der mathematischen Wissenschaften, Bd. {\bf 254}. Springer, Berlin Heidelberg  (1983)
\medskip
\item{[FS1]} Freitag, E. Salvati Manni, R.:
{\it Some Siegel threefolds with a Calabi-Yau model,\/}
Ann.\ Scuola Norm.\ Sup.\ Pisa Cl.\ Sci.\ (5) Vol.~IX,  833--850 (2010)\hfill\break
arXiv: 0905.415
\medskip
\item{[FS2]} Freitag, E. Salvati Manni, R.:
{\it Some Siegel threefolds with a Calabi-Yau model II,\/}
arXiv.org: 1001.0324  (2010)
\medskip
\item{[GN]} van Geemen, B., Nygaard, N.O.:
{\it On the geometry and arithmetic of some Siegel modular threefolds,\/}
Journal of Number Theory {\bf 53}, 45--87 (1995)
\medskip
\item{[Ko]} Kollar, J.:{\it Flops,\/}
Nagoya Math. J. {\bf 113}, 15--36 (1989)
\medskip
\item{[Re]} Reid, M.: {\it La correspondence de McKay,\/}
S\'eminaire Bourbaki 1999/2000. Ast\'erisque No. {\bf 276}, 53--72 (2002)
\medskip
\item{[Ro]} Roan, S.: {\it Minimal resolutions of Gorenstein
orbifolds in dimension three.\/}
Topology {\bf 35}, 489--508 (1996)
\bye